\documentclass[11pt,a4paper,titlepage]{article}
\usepackage[mathscr,mathcal]{eucal}
\usepackage{colortbl}
\usepackage{amsfonts}
\usepackage{amssymb}
\newtheorem{theorem}{Theorem}[section]
\newtheorem{definition}[theorem]{Definition}

\newtheorem{lemma}[theorem]{Lemma}
\newtheorem{corollary}[theorem]{Corollary}
\newtheorem{proposition}[theorem]{Proposition}

\newtheorem{px}[theorem]{Example}

\newcommand{\ie}{\emph{i.e., }}

\begin{document}
\begin{center}
\textbf{DEGREE AND VALUATION OF THE SCHUR ELEMENTS
OF CYCLOTOMIC HECKE ALGEBRAS}\\$ $\\
                                 
MARIA CHLOUVERAKI\end{center} 
$ $\\  
ABSTRACT. Following the generalization of the notion of families of characters, defined by Lusztig for Weyl groups, to the case of complex reflection groups, thanks to the definition given by Rouquier, we show that the degree and the   valuation of the Schur elements (functions $A$ and $a$) remain constant on the ``families'' of the cyclotomic Hecke algebras of the exceptional complex reflection groups. The same result has already been obtained for the groups of the infinite series and for some special cases of exceptional groups.  \\
   \\
  \textbf{ Acknowledgements. }I would like to thank Jean Michel for making my algorithm look better and run faster. I would also like to thank the Ecole Polytechnique F\'{e}d\'{e}rale de Lausanne for its financial support.

\section*{Introduction}

The work of G. Lusztig on the irreducible characters of reductive
groups over finite fields has displayed the important role of the
``families of characters'' of the Weyl groups concerned. More recent results of Gyoja $\cite{Gy}$ and Rouquier $\cite{Rou}$ have made possible the definition of a substitute for
families of characters which can be applied to all complex
reflection groups. Rouquier has shown that the families
of characters of a Weyl group $W$ are exactly the blocks of
characters of the Iwahori-Hecke algebra of $W$ over a suitable
coefficient ring, the ``Rouquier ring''. This definition generalizes without problem to all
cyclotomic Hecke algebras of complex reflection groups. 

Since the families of characters of the Weyl group play an essential role in the
definition of the families of unipotent characters of the
corresponding finite reductive group (cf. \cite{Lu1}), we can hope that
the families of characters of the cyclotomic Hecke algebras play a key
role in the organization of families of unipotent characters more
generally. Moreover, the determination of these families is crucial for the program 
``Spets'' (cf. \cite{BMM2}), whose ambition is to give to complex reflection groups the role of Weyl groups
of as yet mysterious objects.

In the case of the Weyl groups and their usual Hecke algebra, the families of characters can be defined using the existence of Kazhdan-Lusztig bases. Lusztig attaches to every irreducible character two integers, denoted by $a$ and $A$, and shows (cf. \cite{Lu2}, $3.3$ and $3.4$) that they are constant on the families. 
In an analogue way, we can define integers $a$ and $A$ attached to every irreducible character of a cyclotomic Hecke algebra of a complex reflection group. For the groups of the infinite series, it has been shown that $a$ and $A$ are constant on the Rouquier blocks (cf. \cite{BK}, \cite{Ch1}, \cite{Ch2}). Moreover, Malle and Rouquier have proved that $a$ and $A$ are constant on the Rouquier blocks of the ``spetsial'' cyclotomic Hecke algebra of the ``spetsial'' exceptional complex reflection groups (cf. \cite{MaRo}, Thm. 5.1). The aim of this paper is the proof of the same result for all cyclotomic Hecke algebras of all exceptional irreducible complex reflection groups.

In  \cite{Chlou}, we show that
the Rouquier blocks of a cyclotomic Hecke algebra of any complex reflection group $W$ depend on
some numerical data of the group, its ``essential hyperplanes''. These hyperplanes are defined by the 
factorization of the Schur elements of the generic Hecke algebra $\mathcal{H}$  associated to $W$. 
We can associate a partition of the set $\mathrm{Irr}(W)$ of irreducible characters of $W$ to  every essential hyperplane $H$, which we call ``Rouquier blocks associated with the hyperplane $H$''. Following theorem $\ref{main theorem}$ and proposition $\ref{explain AllBlocks}$, these partitions generate the partition of  $\mathrm{Irr}(W)$ into Rouquier blocks. They have been determined for all exceptional irreducible complex reflection groups in \cite{Chlou}. We have stored these data in a computer file and created the GAP function \emph{AllBlocks} which displays them. We have also created  the function \emph{RouquierBlocks} which calculates the Rouquier blocks of a given cyclotomic Hecke algebra.

Let $\phi$ be a cyclotomic specialization and $\mathcal{H}_\phi$ the corresponding cyclotomic Hecke algebra. For every irreducible character, we  define $a$ and $A$ to be, respectively, the valuation and the degree of the corresponding Schur element in  $\mathcal{H}_\phi$. In order to show that $a$ and $A$ are constant on the Rouquier blocks, we introduce the notions of ``generic valuation'' and ``generic degree'' (definition $\ref{generic degree}$). Then corollary $\ref{factor degrees reduced to 0}$ in combination with proposition $\ref{explain AllBlocks}$ imply 
that it is enough to check whether they remain constant on the Rouquier blocks associated with each essential hyperplane. 

We have created a GAP program which verifies that the generic valuation and the generic degree remain constant on the Rouquier blocks for the groups $G_7$, $G_{11}$, $G_{19}$, $G_{26}$, $G_{28}$ and $G_{32}$. We provide the algorithm in section 6.1. Then Clifford theory allows us to extend this result to the groups $G_4, \ldots, G_{22}$ and $G_{25}$. Finally, in section 6.2, we explain why, for  the remaining exceptional irreducible complex reflection groups, it is enough to check whether the functions $a$ and $A$ remain constant on the Rouquier blocks of the ``spetsial'' cyclotomic Hecke algebra. 

All computer algorithms presented in this article require the GAP package CHEVIE, where, together with Jean Michlel, we have programmed the generic Schur elements of  the exceptional irreducible complex reflection groups in a factorized form (functions \emph{SchurModels} and \emph{SchurData}). Moreover, they require the GAP functions \emph{AllBlocks} and \emph{RouquierBlocks} contained in the file ``RouquierBlocks.g". All the above, along with a program implementing the algorithms, can be found on my webpage: \\ \emph{http://www.math.jussieu.fr/$\sim$chlouveraki}. 

\section{Generalities on blocks}

Let  $\mathcal{O}$ be a Noetherian and integrally closed domain with field of fractions $F$. Let $A$ be an $\mathcal{O}$-algebra free and finitely generated as an $\mathcal{O}$-module. 

\begin{definition}
The blocks of $A$ are the central primitive idempotents of $A$.
\end{definition}

Let $K$ be a finite Galois extension of $F$ such that the algebra $KA:=K \otimes_{\mathcal{O}}A$ is split semisimple. Then there exists a bijection between the set $\mathrm{Irr}(KA)$ of irreducible characters of $KA$ and the set $\mathrm{Bl}(KA)$ of blocks of $KA$ which sends every irreducible character $\chi$ to the central primitive idempotent $e_\chi$.

\begin{theorem}\label{minimality of blocks}\
\begin{enumerate}
  \item We have $1=\sum_{\chi \in \mathrm{Irr}(KA)}e_\chi$
   and the set $\{e_\chi\}_{\chi \in \mathrm{Irr}(KA)}$ is the set of all the blocks of the algebra $KA$.
  \item There exists a unique partition $\mathrm{Bl}(A)$ of
  $\mathrm{Irr}(KA)$ such that
  \begin{description}
    \item[(a)] For all $B \in \mathrm{Bl}(A)$, the idempotent
    $e_B:=\sum_{\chi \in B}e_\chi$ is a block of $A$.
    \item[(b)] We have $1=\sum_{B \in \mathrm{Bl}(A)}e_B$ and for
    every central idempotent $e$ of $A$, there exists a subset
    $\mathrm{Bl}(A,e)$ of $\mathrm{Bl}(A)$ such that
    $$e=\sum_{B \in \mathrm{Bl}(A,e)}e_B.$$
     \end{description}
    In particular the set $\{e_B\}_{B \in \mathrm{Bl}(A)}$ is the set of all the blocks of $A$.
   If $\chi \in B$ for some $B \in \mathrm{Bl}(A)$, we say that
  ``$\chi$ belongs to the block $e_B$''.\  
\end{enumerate}
\end{theorem}

Now let us suppose that there exists a symmetrizing form on $A$, \ie a linear map
$t:A \rightarrow \mathcal{O}$ such that
\begin{itemize}
\item $t(aa')=t(a'a)$ for all $a,a' \in A$,
\item the map
     $$\begin{array}{cccc}
     \hat{t}: & A & \rightarrow & \textrm{Hom}_\mathcal{O}(A,\mathcal{O}) \\
              & a & \mapsto & (x \mapsto t(ax))
     \end{array}$$
     is an isomorphism of $A$-modules-$A$.
\end{itemize}
Then we have the following result due to Geck (cf. \cite{Ge}).

\begin{proposition}\label{schur elements and idempotents}\
\begin{enumerate}
  \item We have
  $$t=\sum_{\chi \in \mathrm{Irr}(KA)}\frac{1}{s_\chi}\chi,$$
  where $s_\chi$ is the Schur element associated to $\chi$.
  \item For all $\chi \in \mathrm{Irr}(KA)$, the central primitive
  idempotent associated to $\chi$ is
  $$e_\chi=\frac{\hat{t}^{-1}(\chi)}{s_\chi}.$$
\end{enumerate}
\end{proposition}

\section{Generic Hecke algebras}

Let $\mu_\infty$ be the group of all the roots of unity in
$\mathbb{C}$ and $K$ a number field contained in
$\mathbb{Q}(\mu_\infty)$. We denote by $\mu(K)$ the group of all the
roots of unity of $K$. For every integer $d>1$, we set
$\zeta_d:=\mathrm{exp}(2\pi i/d)$ and denote by $\mu_d$ the group of
all the $d$-th roots of unity. 

Let $V$ be a $K$-vector space of
finite dimension $r$. Let $W$ be a finite subgroup of $\mathrm{GL}(V)$ generated by
(pseudo-)reflections acting irreducibly on $V$. Let us denote by $\mathcal{A}$ the set of the
reflecting hyperplanes of $W$.  For every orbit $\mathcal{C}$ of $W$ on $\mathcal{A}$, we denote by
$e_{\mathcal{C}}$ the common order of the subgroups $W_H$, where $H$
is any element of $\mathcal{C}$ and $W_H$ the subgroup formed by $\mathrm{id}_V$
and all the reflections fixing the hyperplane $H$.

We choose a set of indeterminates
$\textbf{u}=(u_{\mathcal{C},j})_{(\mathcal{C} \in
\mathcal{A}/W)(0\leq j \leq e_{\mathcal{C}}-1)}$ and we denote by
$\mathbb{Z}[\textbf{u},\textbf{u}^{-1}]$ the Laurent polynomial ring
in all the indeterminates $\textbf{u}$. If we denote by $B$ the braid group associated to $W$, then we define the \emph{generic
Hecke algebra} $\mathcal{H}$ of $W$ to be the quotient of the group
algebra $\mathbb{Z}[\textbf{u},\textbf{u}^{-1}]B$ by the ideal
generated by the elements of the form
$$(\textbf{s}-u_{\mathcal{C},0})(\textbf{s}-u_{\mathcal{C},1}) \ldots (\textbf{s}-u_{\mathcal{C},e_{\mathcal{C}}-1}),$$
where $\mathcal{C}$ runs over the set $\mathcal{A}/W$ and
$\textbf{s}$ runs over the set of monodromy generators around the
images in $\mathcal{M}/W$ of the elements of the hyperplane
orbit $\mathcal{C}$.

\begin{px}\label{first ex}
\small{\emph{Let $W:=G_4=<s,t \,|\, \,sts=tst, s^3=t^3=1>$. Then $s$
and $t$ are conjugate in $W$ and their reflecting hyperplanes belong
to the same orbit in $\mathcal{A}/W$. The generic Hecke algebra of
$W$ can be presented as follows
$$\begin{array}{rll}
   \mathcal{H}(G_4)=<S,T \,\,|&STS=TST, &(S-u_0)(S-u_1)(S-u_2)=0, \\
                                &         &(T-u_0)(T-u_1)(T-u_2)=0>.
  \end{array}$$}}
\end{px}

From now on, we assume that the algebra $\mathcal{H}$ is a free
$\mathbb{Z}[\textbf{u},\textbf{u}^{-1}]$-module of
rank $|W|$ and that there exists a symmetrizing form $t$ on $\mathcal{H}$
which satisfies certain conditions (cf., for example, \cite{BK}, Hyp. 2.1). 
Note that
the above assumptions have been verified for all but a finite number of irreducible
complex reflection groups (\cite{BMM2}, remarks before 1.17, $\S$ 2;
\cite{GIM}).

Then
we have the following result by G.Malle (\cite{Ma4}, 5.2).

\begin{theorem}\label{Semisimplicity Malle}
Let $\textbf{\emph{v}}=(v_{\mathcal{C},j})_{(\mathcal{C} \in
\mathcal{A}/W)(0\leq j \leq e_{\mathcal{C}}-1)}$ be a set of
$\sum_{\mathcal{C} \in \mathcal{A}/W}e_{\mathcal{C}}$ indeterminates
such that, for every $\mathcal{C},j$, we have
$$v_{\mathcal{C},j}^{|\mu(K)|}:=\zeta_{e_\mathcal{C}}^{-j}u_{\mathcal{C},j}.$$
The field $K$ is the field of definition of $W$ and the element $\zeta_{e_\mathcal{C}}$ belongs to $K$.
We have that the $K(\textbf{\emph{v}})$-algebra
$K(\textbf{\emph{v}})\mathcal{H}$ is split semisimple.
\end{theorem}

\begin{px}\label{second ex}
\small{\emph{If $W=G_4$ and $K=\mathbb{Q}(\zeta_3)$, then, in the example \ref{first ex},
we replace $u_0,u_1,u_2$ by $v_0^6,\zeta_3v_1^6,\zeta_3^2v_2^6$.
The algebra $\mathbb{Q}(\zeta_3,v_0,v_1,v_2)\mathcal{H}(G_4)$ is split semisimple.}}
\end{px}

By ``Tits' deformation theorem'' (cf., for example, \cite{BMM2}, 7.2), it follows
that the specialization $v_{\mathcal{C},j}\mapsto 1$ induces a
bijection $\chi \mapsto \chi_{\textbf{v}}$ from the
set $\mathrm{Irr}(K(\textbf{v})\mathcal{H})$ of absolutely
irreducible characters of $K(\textbf{v})\mathcal{H}$ to the set
$\mathrm{Irr}(W)$ of absolutely irreducible characters of $W$.

The following result concerning the form of the Schur elements associated
with the irreducible characters of $K(\textbf{v})\mathcal{H}$ is proved in \cite{Chlou}, Thm. 3.2.5, using  case by case analysis (cf. \cite{AlLu}, \cite{Ben}, \cite{Lu82}, \cite{Lu79},  \cite{Ma2}, \cite{Ma5}, \cite{Mat}, \cite{Sur}) and Clifford theory.

\begin{theorem}\label{Schur element generic}
The Schur element $s_\chi(\textbf{\emph{v}})$ associated with the
character $\chi_{\textbf{\emph{v}}}$ of
$K(\textbf{\emph{v}})\mathcal{H}$ is an element of
$\mathbb{Z}_K[\textbf{\emph{v}},\textbf{\emph{v}}^{-1}]$ of the form
$$s_\chi({\textbf{\emph{v}}})=\xi_\chi N_\chi \prod_{i \in I_\chi} \Psi_{\chi,i}(M_{\chi,i})^{n_{\chi,i}}$$
where
\begin{itemize}
    \item $\xi_\chi$ is an element of $\mathbb{Z}_K$,
    \item $N_\chi= \prod_{\mathcal{C},j} v_{\mathcal{C},j}^{b_{\mathcal{C},j}}$ is a monomial in $\mathbb{Z}_K[\textbf{\emph{v}},\textbf{\emph{v}}^{-1}]$
          such that $\sum_{j=0}^{e_\mathcal{C}-1}b_{\mathcal{C},j}=0$
          for all $\mathcal{C} \in \mathcal{A}/W$,
    \item $I_\chi$ is an index set,
    \item $(\Psi_{\chi,i})_{i \in I_\chi}$ is a family of $K$-cyclotomic polynomials in one variable
           (i.e., minimal polynomials of the roots of unity over $K$),
    \item $(M_{\chi,i})_{i \in I_\chi}$ is a family of monomials in $\mathbb{Z}_K[\textbf{\emph{v}},\textbf{\emph{v}}^{-1}]$
          and if $M_{\chi,i} = \prod_{\mathcal{C},j} v_{\mathcal{C},j}^{a_{\mathcal{C},j}}$,
          then $\textrm{\emph{gcd}}(a_{\mathcal{C},j})=1$
          and $\sum_{j=0}^{e_\mathcal{C}-1}a_{\mathcal{C},j}=0$
          for all $\mathcal{C} \in \mathcal{A}/W$,
    \item ($n_{\chi,i})_{i \in I_\chi}$ is a family of positive integers.
\end{itemize}
This factorization is unique in $K[\textbf{\emph{v}},\textbf{\emph{v}}^{-1}]$. Moreover, the monomials
 $(M_{\chi,i})_{i \in I_\chi}$ are unique up to inversion, whereas the coefficient $\xi_\chi$ is unique up to multiplication by a root of unity.
\end{theorem}

\begin{px}\label{third ex}
\small{\emph{Let us denote by $\theta$ the only irreducible character of degree $3$ of $G_4$.
If $v_0,v_1,v_2$ are defined as in example  \ref{second ex}, then we have
$$s_\theta(\textbf{v})=\prod_{i=0}^2\Phi_4(v_i^2 v_{i+1}^{-1}v_{i+2}^{-1})\Phi_{12}'(v_i^2 v_{i+1}^{-1}v_{i+2}^{-1})
\Phi_{12}''(v_i^2 v_{i+1}^{-1}v_{i+2}^{-1}),$$
where $\Phi_{4}(q)=q^2+1$, $\Phi_{12}'(q)=q^2+\zeta_3^2$, $\Phi_{12}''(q)=q^2+\zeta_3$ and the indexes are taken
$\mathrm{mod}\,3$.}}
\end{px}

Let $A:=\mathbb{Z}_K[\textbf{v},\textbf{v}^{-1}]$ and $\mathfrak{p}$ be a prime ideal of $\mathbb{Z}_K$. 

\begin{definition}\label{p-essential monomial}
Let  $M = \prod_{\mathcal{C},j} v_{\mathcal{C},j}^{a_{\mathcal{C},j}}$ be a monomial in $A$
such that $\textrm{\emph{gcd}}(a_{\mathcal{C},j})=1$. We say that $M$ is $\mathfrak{p}$-essential
for a character $\chi \in \mathrm{Irr}(W)$, if there exists a $K$-cyclotomic polynomial $\Psi$ such that
\begin{itemize}
\item $\Psi(M)$ divides $s_\chi(\textbf{\emph{v}})$.
\item $\Psi(1)  \in  \mathfrak{p}$.
\end{itemize}
We say that $M$ is $\mathfrak{p}$-essential 
for $W$, if there exists a character $\chi \in \mathrm{Irr}(W)$ such that
$M$ is $\mathfrak{p}$-essential for $\chi$.
\end{definition}

\begin{px}\label{fourth ex}
\small{\emph{The monomials $v_0^2 v_{1}^{-1}v_{2}^{-1}$, $v_1^2 v_{2}^{-1}v_{0}^{-1}$ and
$v_2^2 v_{0}^{-1}v_{1}^{-1}$ are $2$-essential for the irreducible character 
of degree $3$ of $G_4$.}}
\end{px}

The following proposition (\cite{Chlou}, Prop. 3.2.6) gives a characterization of $\mathfrak{p}$-essential monomials, which plays an essential role in the proof of theorem $\ref{main theorem}$. 

\begin{proposition}\label{p-essential}
Let  $M = \prod_{\mathcal{C},j} v_{\mathcal{C},j}^{a_{\mathcal{C},j}}$ be a monomial in $A$
such that $\textrm{\emph{gcd}}(a_{\mathcal{C},j})=1$. We set  $\mathfrak{q}_M:=(M-1)A +\mathfrak{p}A$.
Then
\begin{enumerate}
\item The ideal $\mathfrak{q}_M$ is a prime ideal of $A$. 
\item $M$ is $\mathfrak{p}$-essential for $\chi \in \mathrm{Irr}(W)$ if and only if
$s_\chi(\textbf{\emph{v}})/\xi_\chi \in \mathfrak{q}_M$.
\end{enumerate}
\end{proposition}

\section{Cyclotomic Hecke algebras}

Let $y$ be an indeterminate. We set $x:=y^{|\mu(K)|}.$

\begin{definition}\label{specialization}
A cyclotomic specialization of $\mathcal{H}$ is a
$\mathbb{Z}_K$-algebra morphism $\phi:
\mathbb{Z}_K[\textbf{\emph{v}},\textbf{\emph{v}}^{-1}]\rightarrow
\mathbb{Z}_K[y,y^{-1}]$ with the following properties:
\begin{itemize}
  \item $\phi: v_{\mathcal{C},j} \mapsto y^{n_{\mathcal{C},j}}$ where
  $n_{\mathcal{C},j} \in \mathbb{Z}$ for all $\mathcal{C}$ and $j$.
  \item For all $\mathcal{C} \in \mathcal{A}/W$, if $z$ is another
  indeterminate, the element of $\mathbb{Z}_K[y,y^{-1},z]$ defined by
  $$\Gamma_\mathcal{C}(y,z):=\prod_{j=0}^{e_\mathcal{C}-1}(z-\zeta_{e_\mathcal{C}}^jy^{n_{\mathcal{C},j}})$$
  is invariant by the action of $\textrm{\emph{Gal}}(K(y)/K(x))$.
\end{itemize}
\end{definition}

If $\phi$ is a cyclotomic specialization of $\mathcal{H}$,
the corresponding \emph{cyclotomic Hecke algebra} is the
$\mathbb{Z}_K[y,y^{-1}]$-algebra, denoted by $\mathcal{H}_\phi$,
which is obtained as the specialization of the
$\mathbb{Z}_K[\textbf{v},\textbf{v}^{-1}]$-algebra $\mathcal{H}$ via
the morphism $\phi$. It also has a symmetrizing form $t_\phi$
defined as the specialization of the canonical form $t$.

\begin{px}\label{spetsial}
\small{\emph{The ``spetsial'' Hecke algebra $\mathcal{H}^s(W)$ is the
cyclotomic algebra obtained by the specialization
$$v_{\mathcal{C},0} \mapsto y,\,\, v_{\mathcal{C},j} \mapsto1 \textrm{
for } 1 \leq j \leq e_\mathcal{C}-1, \textrm{ for all } \mathcal{C}
\in \mathcal{A}/W.$$ }}
\end{px}

The following result is proved in \cite{Chlou} (remarks following Thm. 3.3.3):

\begin{proposition}
The algebra $K(y)\mathcal{H}_\phi$ is split semisimple. 
\end{proposition}

For $y=1$ this algebra specializes to the group
algebra $KW$. Thus, by ``Tits' deformation theorem'', the
specialization $v_{\mathcal{C},j} \mapsto 1$ defines the following bijections
$$\begin{array}{ccccc}
    \textrm{Irr}(K(\textbf{v})\mathcal{H}) & \leftrightarrow & \textrm{Irr}(K(y)\mathcal{H}_\phi) & \leftrightarrow & \textrm{Irr}(W) \\
    \chi_{\textbf{v}} & \mapsto & \chi_{\phi} & \mapsto & \chi.
  \end{array}$$

The following result is an immediate consequence of Theorem
$\ref{Schur element generic}$.

\begin{proposition}\label{Schur element cyclotomic}
The Schur element $s_{\chi_\phi}(y)$ associated with the irreducible
character $\chi_\phi$ of $K(y)\mathcal{H}_\phi$ is a Laurent
polynomial in $y$ of the form
$$s_{\chi_\phi}(y)=\psi_{\chi,\phi} y^{a_{\chi,\phi}} \prod_{\Phi \in
C_{\chi,\phi}}\Phi(y)^{n_{\chi,\Phi}}$$ where $\psi_{\chi,\phi} \in
\mathbb{Z}_K$, $a_{\chi,\phi} \in \mathbb{Z}$, $n_{\chi,\Phi} \in
\mathbb{N}$ and $C_{\chi,\phi}$ is a set of $K$-cyclotomic polynomials.
\end{proposition}

\section{Rouquier blocks}

\begin{definition}\label{Rouquier ring}
We call Rouquier ring of $K$ and denote by $\mathcal{R}$ the
$\mathbb{Z}_K$-subalgebra of $K(y)$
$$\mathcal{R}:=\mathbb{Z}_K[y,y^{-1},(y^n-1)^{-1}_{n\geq 1}]$$
\end{definition}

Let $\phi: v_{\mathcal{C},j} \mapsto y^{n_{\mathcal{C},j}}$ be a
cyclotomic specialization and $\mathcal{H}_\phi$ the corresponding
cyclotomic Hecke algebra. The \emph{Rouquier blocks} of
$\mathcal{H}_\phi$ are the blocks of the algebra
$\mathcal{R}\mathcal{H}_\phi$.\\$ $\\
\begin{remark} \emph{It has been shown by Rouquier \cite{Rou} that if $W$ is a Weyl group and $\mathcal{H}_\phi$ is the ``spetsial'' cyclotomic Hecke algebra (see ex. $\ref{spetsial}$), then its Rouquier blocks coincide with the ``families of characters'' defined by Lusztig.}
\end{remark}\\

Due to the form of the cyclotomic Schur elements, the form of the prime ideals of the Rouquier ring
(see, for example, \cite{Chlou}, Prop. 3.4.2) and an elementary result of blocks theory (see, for example, \cite{BK}, Prop. 1.13), we obtain the following description of the Rouquier blocks:

\begin{proposition}\label{Rouquier blocks and central characters}
Two characters $\chi,\psi \in \emph{Irr}(W)$ are in the same Rouquier block of $\mathcal{H}_\phi$
if and only if there exists a finite sequence
$\chi_0,\chi_1,\ldots,\chi_n \in \emph{Irr}(W)$ and a finite
sequence $\mathfrak{p}_1,\ldots,\mathfrak{p}_n$ of prime
ideals of $\mathbb{Z}_K$ such that
\begin{itemize}
  \item $\chi_0=\chi$ and $\chi_n=\psi$,
  \item for all $j$ $(1\leq j \leq n)$,\,\,
         the characters $\chi_{j-1}$ and $\chi_{j}$
         belong to the same block of
         $\mathcal{R}_{\mathfrak{p}_j\mathcal{R}}\mathcal{H}_\phi.$
 \end{itemize}                   
\end{proposition}

The above proposition implies that if we know the blocks of the algebra $\mathcal{R}_{\mathfrak{p}\mathcal{R}}\mathcal{H}_\phi$ for every  prime ideal
$\mathfrak{p}$  of $\mathbb{Z}_K$,
then we know the Rouquier blocks of  $\mathcal{H}_\phi$. In order to determine the former, we
can use the following theorem (\cite{Cras}, Thm. 2.5)

\begin{theorem}\label{main theorem}
Let $A:=\mathbb{Z}_K[\textbf{\emph{v}},\textbf{\emph{v}}^{-1}]$ and $\mathfrak{p}$ be a prime
ideal of $\mathbb{Z}_K$. 
Let $M_1,\ldots,M_k$ be all the
$\mathfrak{p}$-essential monomials for $W$ such that $\phi(M_j)=1$
for all $j=1,\ldots,k$. Set $\mathfrak{q}_0:=\mathfrak{p}A$,
$\mathfrak{q}_j:=\mathfrak{p}A+(M_j-1)A$ for $j=1,\ldots,k$ and
$\mathcal{Q}:=\{\mathfrak{q}_0,\mathfrak{q}_1,\ldots,\mathfrak{q}_k\}$.
Two irreducible characters $\chi,\psi \in \textrm{\emph{Irr}}(W)$
are in the same block of $\mathcal{R}_{\mathfrak{p}\mathcal{R}}\mathcal{H}_\varphi$ if
and only if there exist a finite sequence
$\chi_0,\chi_1,\ldots,\chi_n \in \textrm{\emph{Irr}}(W)$ and a
finite sequence $\mathfrak{q}_{j_1},\ldots,\mathfrak{q}_{j_n} \in
\mathcal{Q}$ such that
\begin{itemize}
  \item $\chi_0=\chi$ and $\chi_n=\psi$,
  \item for all $i$ $(1\leq i \leq n)$,  the characters $\chi_{i-1}$ and $\chi_i$ are
  in the same block of $A_{\mathfrak{q}_{j_i}}\mathcal{H}$.
\end{itemize}
\end{theorem}

Let $\mathfrak{p}$ be a prime ideal of $\mathbb{Z}_K$ and  $\phi: v_{\mathcal{C},j} \mapsto y^{n_{\mathcal{C},j}}$ a cyclotomic specialization.
If $M=\prod_{\mathcal{C},j}v_{\mathcal{C},j}^{a_{\mathcal{C},j}}$
is a $\mathfrak{p}$-essential monomial for $W$, then
$$\phi(M)=1 \Leftrightarrow \sum_{\mathcal{C},j}a_{\mathcal{C},j}n_{\mathcal{C},j}=0.$$
The
hyperplane defined in $\mathbb{C}^{\sum_{\mathcal{C}\in \mathcal{A}/W}e_\mathcal{C}}$  by the relation
$$\sum_{\mathcal{C},j}a_{\mathcal{C},j}t_{\mathcal{C},j}=0,$$ where
$(t_ {\mathcal{C},j})_{ \mathcal{C},j}$ is a set of $\sum_{\mathcal{C}\in \mathcal{A}/W}e_\mathcal{C}$
indeterminates, is called \emph{$\mathfrak{p}$-essential hyperplane}
for $W$. A hyperplane in $\mathbb{C}^{\sum_{\mathcal{C}\in \mathcal{A}/W}e_\mathcal{C}}$ is called \emph{essential}
for $W$, if it is $\mathfrak{p}$-essential for some prime ideal
$\mathfrak{p}$ of $\mathbb{Z}_K$.

\begin{px}\label{fifth ex}
\small{\emph{The essential hyperplanes of $G_4$ are:
$$H_{0,1} : t_0-t_1=0, \,  \,H_{0,2} : t_0-t_2=0,  \, \, H_{1,2} : t_1-t_2=0,$$
$$H_0 : 2t_0-t_1-t_2=0, \, \, H_1 : 2t_1-t_2-t_0=0,  \, \,H_2 : 2t_2-t_0-t_1=0.$$}}
\end{px}

Let $\phi_\emptyset: v_{\mathcal{C},j} \mapsto y^{n_{\mathcal{C},j}}$ be a cyclotomic specialization such that the $n_{\mathcal{C},j}$ belong to no essential hyperplane. We call \emph{Rouquier blocks associated with no essential hyperplane} and denote by $\mathcal{B}^\emptyset$ the partition of $\mathrm{Irr}(W)$ into the Rouquier blocks of $\mathcal{H}_{\phi_\emptyset}$. Now let $H$ be an essential hyperplane for $W$ and let $\phi_H: v_{\mathcal{C},j} \mapsto y^{n_{\mathcal{C},j}}$ be a cyclotomic specialization such that the $n_{\mathcal{C},j}$ belong to the essential hyperplane $H$ and no other.  We call \emph{Rouquier blocks associated with the essential hyperplane $H$} and denote by $\mathcal{B}^H$ the partition of $\mathrm{Irr}(W)$ into the Rouquier blocks of $\mathcal{H}_{\phi_H}$. Due to theorem $\ref{main theorem}$, the partition $\mathcal{B}^H$ is coarser than the partition $\mathcal{B}^\emptyset$.

The following result is an immediate consequence of 
proposition $\ref{Rouquier blocks and central characters}$ and theorem $\ref{main theorem}$.

\begin{proposition}\label{explain AllBlocks}
Let  $\phi: v_{\mathcal{C},j} \mapsto y^{n_{\mathcal{C},j}}$ be a cyclotomic specialization. 
 If the $n_{\mathcal{C},j}$ belong to no essential hyperplane for $W$, then 
the Rouquier blocks of $\mathcal{H}_\phi$ coincide with the partition  $\mathcal{B}^\emptyset$.
Otherwise, let $\mathcal{E}$ be the set of all essential hyperplanes that the   $n_{\mathcal{C},j}$ belong to. Two irreducible characters $\chi, \psi \in \mathrm{Irr}(W)$ belong to the same Rouquier block of 
$\mathcal{H}_\phi$ if and only if there exist 
a finite sequence
$\chi_0,\chi_1,\ldots,\chi_n \in \textrm{\emph{Irr}}(W)$ and a
finite sequence $H_{1},\ldots,H_{n} \in
\mathcal{E}$ such that
\begin{itemize}
  \item $\chi_0=\chi$ and $\chi_n=\psi$,
  \item for all $i$ $(1\leq i \leq n)$,  the characters $\chi_{i-1}$ and $\chi_i$ belong to the same part of
  $\mathcal{B}^{H_i}$.
\end{itemize}
\end{proposition}

\begin{px}\label{sixth ex}
\small{\emph{Let $\phi^s : v_i \mapsto y^{n_i}$ be the ``spetsial'' cyclotomic specialization for $G_4$, \ie
$n_0=1$ and $n_1=n_2=0$. The integers $n_i$ belong only to the essential hyperplane
$H_{1,2} : t_1-t_2=0$ and therefore, the Rouquier blocks of $\mathcal{H}^s$ coincide
with the partition $\mathcal{B}^{H_{1,2}}$.}}
\end{px}

In the fourth chapter of \cite{Chlou}, we explain how we have obtained the partitions $\mathcal{B}^\emptyset$ and $\mathcal{B}^H$ for every essential hyperplane $H$ for every exceptional irreducible complex reflection group. 

\section{Functions $a$ and $A$}

Following the notations in \cite{BMM2}, 6B, for every element $P(y)
\in \mathbb{C}(y)$, we call
\begin{itemize}
  \item \emph{valuation of $P(y)$ at $y$} and denote by $\mathrm{val}_y(P)$ the order of $P(y)$
  at 0 (we have $\mathrm{val}_y(P)<0$ if 0 is a pole of $P(y)$ and $\mathrm{val}_y(P)>0$ if 0 is a zero of $P(y)$),
  \item \emph{degree of $P(y)$ at $y$} and denote by $\mathrm{deg}_y(P)$ the negative of the
  valuation of $P(1/y)$.
\end{itemize}
For $\chi
\in \mathrm{Irr}(W)$, we define
$$a_{\chi_\phi}:=\mathrm{val}_y(s_{\chi_\phi}(y)) \,\textrm{ and }\,
A_{\chi_\phi}:=\mathrm{deg}_y(s_{\chi_\phi}(y)).$$ 
The proof of the following
result can be found in \cite{BK}, Prop. 2.9.

\begin{proposition}\label{aA}\
Let $\chi,\psi \in \mathrm{Irr}(W)$. If $\chi_\phi$ and
        $\psi_\phi$ belong to the same Rouquier block, then
        $$a_{\chi_\phi}+A_{\chi_\phi}=a_{\psi_\phi}+A_{\psi_\phi}.$$
\end{proposition}

In the next section, we are going to prove that the functions $a$ and $A$ are constant on the Rouquier blocks of the cyclotomic Hecke algebras of the exceptional irreducible complex reflection groups. In order to do that, we need to prove some results concerning the valuation and the degree of the Schur elements which hold for all complex reflection groups. First, let us introduce the symbols  $(y^n)^+$
and  $(y^n)^-$.

\begin{definition}\label{symbols} Let  $n \in \mathbb{Z}$.
\begin{itemize}
  \item $(y^n)^+ = \left\{
                      \begin{array}{ll}
                        n, & \hbox{if $n > 0$;} \\
                        0, & \hbox{if $n \leq 0$.}
                      \end{array}
                    \right.$
  \item $(y^n)^- = \left\{
                      \begin{array}{ll}
                        n, & \hbox{if $n < 0$;} \\
                        0, & \hbox{if $n \geq 0$.}
                      \end{array}
                    \right.$
\end{itemize}
\end{definition}

Now let  $\chi \in
\mathrm{Irr}(W)$. Following the notations of Theorem $\ref{Schur element generic}$, the generic Schur 
element $s_\chi(\textbf{v})$
associated to $\chi$ is an element of
$\mathbb{Z}_K[\textbf{v},\textbf{v}^{-1}]$ of the form
$$s_\chi(\textbf{v})=\xi_\chi N_\chi \prod_{i \in I_\chi} \Psi_{\chi,i}(M_{\chi,i})^{n_{\chi,i}}.\,\,\,\,\,\,\,\,\,\,\,\,\,\,\,\,(\dag)$$
We fix the factorization $(\dag)$ for $s_\chi(\textbf{v}).$

\begin{proposition}\label{Aa formula}
Let $\phi:v_{\mathcal{C},j} \mapsto y^{n_{\mathcal{C},j}}$ be a
cyclotomic specialization. Then
\begin{itemize}
  \item $a_{\chi_\phi}=
  \sum_{\mathcal{C},j}b_{\mathcal{C},j}n_{\mathcal{C},j}+
  \sum_{i \in I_\chi}
  n_{\chi,i}\mathrm{deg}(\Psi_{\chi,i})(\phi(M_{\chi,i}))^-$.
  \item $A_{\chi_\phi}=\sum_{\mathcal{C},j}b_{\mathcal{C},j}n_{\mathcal{C},j}+
  \sum_{i \in I_\chi}
  n_{\chi,i}\mathrm{deg}(\Psi_{\chi,i})(\phi(M_{\chi,i}))^+$.
\end{itemize}
\end{proposition}

\begin{px}\label{seventh ex}
\small{\emph{Let $\phi$ be a cyclotomic specialization for $G_4$.
Following the factorization of the generic Schur element of the character $\theta$ in example $\ref{third ex}$, we have that
\begin{itemize}
  \item $a_{\theta_\phi}=6\cdot (\phi(v_0^2 v_{1}^{-1}v_{2}^{-1})^-+
  \phi(v_1^2 v_{2}^{-1}v_{0}^{-1})^-+\phi(v_2^2 v_{0}^{-1}v_{1}^{-1})^-)$,
  \item $A_{\theta_\phi}=6\cdot (\phi(v_0^2 v_{1}^{-1}v_{2}^{-1})^++
  \phi(v_1^2 v_{2}^{-1}v_{0}^{-1})^++\phi(v_2^2 v_{0}^{-1}v_{1}^{-1})^+)$.
\end{itemize}
If $\phi$ is the ``spetsial'' cyclotomic specialization, \ie $\phi(v_0)=y$ and $\phi(v_1)=\phi(v_2)=1$, then
$$\begin{array}{lll}
\phi(v_0^2 v_{1}^{-1}v_{2}^{-1})^-=0, &\phi(v_1^2 v_{2}^{-1}v_{0}^{-1})^-= -1,
&\phi(v_2^2 v_{0}^{-1}v_{1}^{-1})^-=-1,\\
\phi(v_0^2 v_{1}^{-1}v_{2}^{-1})^+=2, &\phi(v_1^2 v_{2}^{-1}v_{0}^{-1})^+= 0,
&\phi(v_2^2 v_{0}^{-1}v_{1}^{-1})^+=0,
\end{array}$$
Thus, we have $a_{\theta_\phi}=-12$ and  $A_{\theta_\phi}=12$.}}
\end{px}

\begin{definition}\label{factor degree}
Let $M=\prod_{\mathcal{C},j} v_{\mathcal{C},j}^{a_{\mathcal{C},j}}$ be a monomial 
 with $\textrm{\emph{gcd}}(a_{\mathcal{C},j})=1$
 and $\Psi$ a $K$-cyclotomic polynomial such that $\Psi(M)$ appears in $(\dag)$. 
 The factor degree of $\Psi(M)$ for $\chi$ with respect to $(\dag)$ is defined as the product
 $$f_{\Psi(M)}(\textbf{\emph{t}})= \mathrm{deg}(\Psi)\cdot(\sum_{\mathcal{C},j} a_{\mathcal{C},j}{t_{\mathcal{C},j}}),$$
 where $\textbf{\emph{t}}=(t_{\mathcal{C},j})_{\mathcal{C},j}$ is a set of $\sum_{\mathcal{C} \in \mathcal{A}/W}e_\mathcal{C}$ indeterminates. If $n$ is the greatest positive integer such that $\Psi(M)^n$ appears in $(\dag)$, then $n$ is called the coefficient of the factor degree $f_{\Psi(M)}$ and it is denoted by $\textbf{\emph{c}}(f_{\Psi(M)})$.
 \end{definition}
 
 Then we can define an equivalence relation on the set $\mathcal{F}_\chi $ of all factor degrees for $\chi$  with respect to $(\dag)$: 
 
 \begin{definition}\label{equivalence}
 Two factor degrees $f_1,f_2$ are equivalent, if there exists a positive number $q \in \mathbb{Q}$ such that $f_1=qf_2$. We write $f_1 \sim f_2$.
 \end{definition}

 \begin{definition}\label{good sign map}
 Let $\mathcal{F}_\chi$ be the set of all factor degrees for  $\chi$ with respect to $(\dag)$ and let $\epsilon$ be a sign map for $\mathcal{F}_\chi$, i.e., a map
 $\mathcal{F}_\chi \rightarrow \{-1,1\}$. We say that $\epsilon$ is a good sign map for $\mathcal{F}_\chi$ if it satisfies the following conditions:
\begin{enumerate}
\item If $f_1,f_2 \in \mathcal{F}_\chi$ with $f_1 \sim f_2$, then $\epsilon(f_1)=\epsilon(f_2).$
\item If $f_1,f_2 \in \mathcal{F}_\chi$ with $f_1 \sim -f_2$, then $\epsilon(f_1)=-\epsilon(f_2).$
\end{enumerate}
\end{definition}
 
 In order to obtain the main result, we need to introduce the notions of  generic valuation and  generic degree of the Schur element   $s_\chi(\textbf{v})$. 
 
 \begin{definition}\label{generic degree}
 Let $\mathcal{F}_\chi$ be the set of all factor degrees for  $\chi$ with respect to $(\dag)$ and let
 $\epsilon: \mathcal{F}_\chi \rightarrow \{-1,1\}$ be a good sign map for $\mathcal{F}_\chi$. Then 
 \begin{itemize}
 \item the generic valuation
 $a_{\chi,\epsilon}(\textbf{\emph{t}})$ of
 $s_\chi(\textbf{\emph{v}})$ with respect  to $\epsilon$ is
 $$a_{\chi,\epsilon}(\textbf{\emph{t}}) :=  \sum_{\mathcal{C},j}b_{\mathcal{C},j}t_{\mathcal{C},j}+
      \sum_{\{f \in \mathcal{F}_\chi \,|\, \epsilon(f)=-1\}}\textbf{\emph{c}}(f) \cdot f.$$
 \item the generic degree 
 $A_{\chi,\epsilon}(\textbf{\emph{t}})$ of
 $s_\chi(\textbf{\emph{v}})$ with respect to $\epsilon$ is
 $$A_{\chi,\epsilon}(\textbf{\emph{t}}):=  \sum_{\mathcal{C},j}b_{\mathcal{C},j}t_{\mathcal{C},j}+
      \sum_{\{f \in \mathcal{F}_\chi \,|\, \epsilon(f)=1\}}\textbf{\emph{c}}(f) \cdot f.$$
 \end{itemize}
 \end{definition}     

The following result is a consequence of the  above definitions  and proposition $\ref{Aa formula}$.
      
 \begin{proposition}\label{generic cyclotomic}
 Let $\phi:v_{\mathcal{C},j} \mapsto y^{n_{\mathcal{C},j}}$ be a cyclotomic specialization and $\chi,\psi \in \mathrm{Irr}(W)$ with sets of factor degrees $\mathcal{F}_\chi, \mathcal{F}_\psi$ respectively.  
 If $a_{\chi,\epsilon}(\textbf{\emph{t}}) = a_{\psi,\epsilon}(\textbf{\emph{t}})$ (resp. $A_{\chi,\epsilon}(\textbf{\emph{t}})= 
 A_{\psi,\epsilon}(\textbf{\emph{t}})$) with respect to every good sign map $\epsilon$ for $\mathcal{F}_\chi \cup \mathcal{F}_\psi$, then $a_{\chi_\phi}=a_{\psi_\phi}$ (resp. $A_{\chi_\phi}=A_{\psi_\phi})$ .
 \end{proposition}
 \begin{apod}{Let $\textbf{n}:=(n_{\mathcal{C},j})_{\mathcal{C},j}$. 
 There exists a good sign map $\epsilon$ for $\mathcal{F}_\chi \cup \mathcal{F}_\psi$ such that
 $$\epsilon(f)=-1 \Leftrightarrow f(\textbf{n}) \leq 0.$$
Then, by proposition  $\ref{Aa formula}$, we have that 
 $$a_{\chi_\phi}=a_{\chi,\epsilon}(\textbf{n})=
 a_{\psi,\epsilon}(\textbf{n})=a_{\psi_\phi}$$ 
 and
 $$A_{\chi_\phi}=A_{\chi,\epsilon}(\textbf{n})=
 A_{\psi,\epsilon}(\textbf{n})=A_{\psi_\phi}.$$ }
 \end{apod}
 
 \begin{corollary}\label{factor degrees reduced to 0}
Let $\phi:v_{\mathcal{C},j} \mapsto y^{n_{\mathcal{C},j}}$ be a cyclotomic specialization such that the integers $n_{\mathcal{C},j}$ belong to the essential hyperplane $H : \sum_{\mathcal{C},j} a_{\mathcal{C},j}{n_{\mathcal{C},j}}=0$. Then we can assume that the set $\textbf{\emph{t}}$ is not algebraically independent, but satisfies $\sum_{\mathcal{C},j} a_{\mathcal{C},j}{t_{\mathcal{C},j}}=0$. If $a_{\chi,\epsilon}(\textbf{\emph{t}}) = a_{\psi,\epsilon}(\textbf{\emph{t}})$ (resp. $A_{\chi,\epsilon}(\textbf{\emph{t}})= 
 A_{\psi,\epsilon}(\textbf{\emph{t}})$) with respect to every good sign map $\epsilon$ for $\mathcal{F}_\chi \cup \mathcal{F}_\psi$, then $a_{\chi_\phi}=a_{\psi_\phi}$ (resp. $A_{\chi_\phi}=A_{\psi_\phi})$.
 \end{corollary}
 
 \section{Exceptional complex reflection groups}

In this section we will prove the following result 
 
\begin{theorem}\label{invariance}
Let $W$ be an exceptional irreducible complex reflection group.
Let $\phi:v_{\mathcal{C},j} \mapsto y^{n_{\mathcal{C},j}}$ be a cyclotomic specialization and $\chi,\psi \in \mathrm{Irr}(W)$. If $\chi_\phi$ and
        $\psi_\phi$ belong to the same Rouquier block, then
        $$a_{\chi_\phi}=a_{\psi_\phi} \textrm{ and } A_{\chi_\phi}=A_{\psi_\phi}.$$
\end{theorem} 
 
\subsection{The groups $G_4, \ldots, G_{22}$, $G_{25}$, $G_{26}$, $G_{28}$, $G_{32}$}
 
 Let $W:=G_m$, where $m \in \{7,11,19,26,28,32\}$. We have created the following algorithm which verifies that the assumptions of corollary $\ref{factor degrees reduced to 0}$ are satisfied on the Rouquier blocks associated with each essential hyperplane.
This algorithm requires the GAP package CHEVIE  and the function $AllBlocks$ contained in the file ``RouquierBlocks.g''.  A program implementing this algorithm can be found on my webpage.\\
\\
\textbf{Algorithm}
\begin{enumerate}
\item We assume that there exists a function $ismultiple(g,f)$ which takes two  polynomials $f,g$ and returns
\begin{itemize}
\item $1$, if there exists a rational $q>0$ such that $g=q*f$,
\item $-1$, if there exists a rational $q<0$ such that $g=q*f$,
\item  $0$, otherwise.
\end{itemize}

\item We define a function $FactorDegrees(H,\chi)$, where
\begin{itemize}
\item $H$ is either the list of the coefficients $a_{\mathcal{C},j}$ of the indeterminates in an essential hyperplane for $W$ or the empty list in the case of  ``no essential hyperplane'',
\item $\chi \in \mathrm{Irr}(W)$  is represented by its position in the list of characters of $W$.
\end{itemize}
In the GAP package CHEVIE,  the functions \emph{SchurModels} and \emph{SchurData} provide us with the irreducible factors and the coefficients of the generic Schur elements of $W$. 
The function $FactorDegrees(H,\chi)$ returns a pair $[F,C]$, where
$F$ is the list of factor degrees of the Schur element of $\chi$ (a list of polynomials)
and $C$ is the term of the generic valuation (and generic degree) induced by the monomial factor $N_\chi$.

\item We assume that there exists a function \emph{SymmetricDifferenceWithMultiplicities}$(l_1,l_2),$ where $l_1$, $l_2$ are two lists, which returns a sublist  $l$ of $l_1 \cup l_2$ such that: $x \in l$ if and only if the multiplicity of $x$ in $l_1$ is different than the multiplicity of $x$ in $l_2$. 

\item
The function $compare(a,b)$  will check the assumptions of corollary $\ref{factor degrees reduced to 0}$ for two irreducible characters $\chi,\psi$. It returns ``true'', if they are satisfied. In order to do that, it takes the corresponding outputs of the function $FactorDegrees$, $a:=[F_\chi,C_\chi]$ and  $b:=[F_\psi,C_\psi]$, and sets $l:=$\emph{SymmetricDifferenceWithMultiplicities}$(F_\chi,F_\psi)$.

If $l$ is empty, then the function returns ``true''. If not, then we have to generate all good sign maps only for $l$, since the common factors don't affect the result:\\ \\
\textbf{Step 1: }We create a sublist $k$ of $l$ such that:
\begin{enumerate}
\item every element of $l$ is a multiple by a non-zero rational number of an element of $k$,
\item if $f,g \in k$ then \emph{ismultiple$(f,g)=0$}.
\end{enumerate}
\textbf{Step 2: }We create a list $a_1$ as follows: For all $f \in F_\chi$, we set $p:=$the position of the $g$ in $k$ such that \emph{ismultiple$(f.g) \neq 0$}. If $p \neq \textrm{false}$, then we add to $a_1$ the triplet
$[f,p,ismultiple(f,k[p])]$. We create a similar list $b_1$ for $F_\psi$.\\ \\  
\textbf{Step 3: }We create all good sign maps for $l$ which consists of creating all the lists of signs of the same length as $k$. Let $M$ be such a matrix and $f \in l$. Then there exists a triplet of the form
$[f,p,ismultiple(k[p],f)]$ in $a_1$ or $b_1$. The corresponding good sign map $\epsilon$ is given by $\epsilon(f):=ismultiple(f,k[p]) \cdot M[p]$. \\ \\
\textbf{Step 4: }We compare 
 $a_{\chi,\epsilon}(\textbf{t})$ with $a_{\psi,\epsilon}(\textbf{t})$ and $A_{\chi,\epsilon}(\textbf{t})$ with $A_{\psi,\epsilon}(\textbf{t})$ (considering only the non-common terms)
 with respect to every good sign map $\epsilon$ for $l$, given that the condition $\sum_{\mathcal{C},j} a_{\mathcal{C},j}{t_{\mathcal{C},j}}=0$ is satisfied.

\item
We create a function $compareblock(H,B)$, where $H$ is a list representing one or no essential hyperplane as in $FactorDegrees$ and $B$ is a block represented as a list of integers, each of which is the position of a character in the list of characters of $W$. If $Length(B)=1$, then it returns '``true''. If not, then it applies $FactorDegrees(H,\chi)$ to all the elements $\chi$ of $B$, creating thus the list \emph{Sch}, and then returns ``true'' if $compare(Sch[1],Sch[j])=true$ for all $j \in \{2,\ldots,Length(B)\}.$

\item
Finally, we create a function $CheckTheorem(m)$ which generates the group $G_m$ and applies
$compareblock(H,B)$ to every $B \in \mathcal{B}^H$, where $H$ runs over the set $\{\emptyset, \textrm{ essential hyperplanes for }W\}$.

\end{enumerate}
The function $CheckTheorem(m)$ has returned ``true'' for all $m \in \{7$, $11$, $19$, $26$, $28$, $32\}$. Then corollary $\ref{factor degrees reduced to 0}$ in combination with proposition $\ref{explain AllBlocks}$, imply that the assertion of Theorem $\ref{invariance}$ holds for $W$.\\

Now let $W:=G_m$,  $m \in \{4,5,6,8,9,10,12,13,14,15,16,17,18,21,22,25\}$. The
fact that Theorem $\ref{invariance}$ holds for the groups $G_7$, $G_{11}$, $G_{19}$ and $G_{26}$ 
and  the use of Clifford theory for the determination of the Schur elements and the Rouquier blocks of the cyclotomic Hecke algebras associated to $W$ (see Appendix)
imply that the assertion of Theorem $\ref{invariance}$ holds for $W$.

\subsection{The other exceptional groups}

Let $W$ be one of the remaining exceptional irreducible complex reflection groups: $G_{23}$, $G_{24}$, $G_{27}$, $G_{29}$,
 $G_{30}$, $G_{31}$, $G_{33}$, $G_{34}$, $G_{35}$, $G_{36}$, $G_{37}$. Then $W$ is generated by reflections of order $2$ whose reflecting hyperplanes belong to one single orbit under the action of $W$. Its generic Hecke algebra is defined over a Laurent polynomial ring in two indeterminates, $v_0$ and $v_1$, and the only possible essential monomial is $v_0v_1^{-1}$. 
Therefore, its generic Schur elements can be expressed as products of $K$-cyclotomic polynomials
in the one variable $v:=v_0v_1^{-1}$, \ie the generic Schur element $s_\chi(v)$
associated to the irreducible character $\chi$  is an element of
$\mathbb{Z}_K[v,v^{-1}]$ of the form
$$s_\chi(v)=\xi_\chi v^{b_\chi} \prod_{\Psi \in C_\chi} \Psi_{\chi}(v)^{n_{\chi,\Psi}},$$
where $\xi_\chi \in \mathbb{Z}_K$, $b_\chi \in \mathbb{Z}$, $C_\chi$ is a set of $K$-cyclotomic polynomials and $n_{\chi,\Psi} \in \mathbb{N}$. If $\phi:v  \mapsto y^n$ ($n \in \mathbb{Z}$) is a
cyclotomic specialization, then
\begin{itemize} 
 \item $a_{\chi_\phi}=n \cdot \mathrm{val}_v(s_\chi(v))$.
 \item $A_{\chi_\phi}=n \cdot \mathrm{deg}_v(s_\chi(v))$.
 \end{itemize}
Therefore, in order to verify theorem $\ref{invariance}$ for $W$, it suffices to check whether the degree and the valuation of the generic Schur elements remain constant on the Rouquier blocks associated with no essential hyperplane. Note that the generic Schur elements coincide with the Schur elements of the ``spetsial'' cyclotomic Hecke algebra
and the Rouquier blocks associated with no essential hyperplane coincide with its Rouquier blocks.

We can easily create an algorithm which returns ``true'' if the degree and the valuation of the Schur elements of the
``spetsial'' cyclotomic Hecke algebra remain constant on its Rouquier blocks. A program realizing this algorithm can be found on my webpage.
It  requires 
the GAP package CHEVIE and the  function $RouquierBlocks$ contained in the file ``RouquierBlocks.g''. Since this program has returned ``true'' for all $m \in \{23,24,27,29,30,31,33,34,35,36,37\}$, we deduce that
the assertion of Theorem $\ref{invariance}$ holds for $W$.

\section{Appendix}

Let us assume that $\mathcal{O}$, $A$ and $K$ are defined as  in section $1$ and that there exists a symmetrizing form $t$ on $A$. 

\begin{definition}\label{symmetric subalgebra}
Let $\bar{A}$ be a subalgebra of $A$ free and of finite
rank as an $\mathcal{O}$-module. We say that $\bar{A}$ is
a symmetric subalgebra of $A$, if it satisfies the following two
conditions:
\begin{enumerate}
  \item $\bar{A}$ is free (of finite rank) as an $\mathcal{O}$-module and the
  restriction $\mathrm{Res}_{\bar{A}}^A(t)$ of the form $t$ to $\bar{A}$ is a symmetrizing form
  for $\bar{A}$,
  \item $A$ is free (of finite rank) as an $\bar{A}$-module for the action
  of left multiplication by the elements of $\bar{A}$.
\end{enumerate}
\end{definition}

From now on, let us suppose that $\bar{A}$ is a symmetric subalgebra
of $A$. Moreover, let $K$ be a finite Galois extension $F$ such that the algebras $KA$ and
$K\bar{A}$ are both split semisimple.

\begin{definition}\label{symmetric algebra of a finite group}
We say that a symmetric $\mathcal{O}$-algebra $(A,t)$ is the twisted
symmetric algebra of a finite group $G$ over the subalgebra
$\bar{A}$, if the following conditions are satisfied:
\begin{itemize}
  \item $\bar{A}$ is a symmetric subalgebra of $A$,
  \item There exists a family $\{A_g \,|\, g \in G\}$ of
  $\mathcal{O}$-submodules of $A$ such that
  \begin{description}
    \item[(a)] $A= \bigoplus_{g \in G}A_g$,
    \item[(b)] $A_1=\bar{A}$,
    \item[(c)] $A_gA_h=A_{gh}$ for all $g,h \in G$,
    \item[(d)] $t(A_g)=0$ for all $g \in G, g \neq  1$,
    \item[(e)] $A_g \cap A^\times \neq \emptyset$ for all $g \in G$ (where
    $A^\times$ is the set of units of $A$).
    \end{description}
\end{itemize}
\end{definition}

\begin{lemma}\label{factor group}
Let $a_g \in A_g$ such that $a_g$ is a unit in $A$. Then
$$A_g=a_g\bar{A}=\bar{A}a_g.$$
\end{lemma}
\begin{apod}{Since $a_g \in A_g$, property (b)
implies that $a_g^{-1} \in A_{g^{-1}}$. If $a \in A_g$, then $
a_g^{-1}a \in A_1=\bar{A}$. We have $a = a_g{a_g}^{-1}a \in
a_g\bar{A}$ and thus $A_g \subseteq a_g\bar{A}$. Property (b)
implies the inverse inclusion. In the same way, we show that
$A_g=\bar{A}a_g$.}
\end{apod}
 
Sometimes the Hecke algebra of a group $W'$ appears as a symmetric
subalgebra of the Hecke algebra of another group $W$, which
contains $W'$. Therefore, it would be helpful,
if we could obtain the Schur elements (resp. the blocks) of the former from the
Schur elements (resp. the blocks) of the
latter. This is possible with the use of a generalization of some
classic results, known as ``Clifford theory'' (see, for example,
\cite{Da}), to the twisted symmetric algebras of finite groups and
more precisely of finite cyclic groups. \\
 
Let $W$ be a complex reflection group and let us denote by
$\mathcal{H}(W)$ its generic Hecke algebra.  Let $\mathcal{H}(W)_\mathrm{sp}$ be 
the algebra obtained from $\mathcal{H}(W)$ by specializing some of the parameters.
Let $W'$ be another complex
reflection group such that $\mathcal{H}(W)_\mathrm{sp}$ is the twisted symmetric algebra of a
finite \textbf{cyclic} group $G$ over the symmetric subalgebra
$\mathcal{H}(W')$. 
In all the cases that will be studied below, applying ``Clifford theory'' (cf., for example,  \cite{BK}, Prop. $1.42$ and $1.45$) gives that
\begin{enumerate}
\item if we denote
by $\chi'$ the (irreducible) restriction to $\mathcal{H}(W')$ of an
irreducible character $\chi \in \mathrm{Irr}(\mathcal{H}(W)_\mathrm{sp})$, then
their Schur elements verify
$$s_\chi = |W:W'| s_{\chi'},$$
\item the blocks of the algebras $\mathcal{H}(W)_\mathrm{sp}$  and $\mathcal{H}(W')$  
          coincide. 
\end{enumerate}
 
\subsection*{The groups $G_4$, $G_5$, $G_6$, $G_7$}

The following table gives the specializations of the parameters of
the generic Hecke algebra $\mathcal{H}(G_7)$,
$(x_0,x_1;y_0,y_1,y_2;z_0,z_1,z_2)$, which give the generic Hecke
algebras of the groups $G_4$, $G_5$ and $G_6$ (\cite{Ma2}, Table
4.6).

\begin{center}
\begin{tabular}{|c|c|c|c|c|}
  \hline
  Group & Index & S & T & U \\
  \hline
  $G_7$ & 1 & $x_0,x_1$ & $y_0,y_1,y_2$         & $z_0,z_1,z_2$\\
  $G_5$ & 2 & $1,-1$    & $y_0,y_1,y_2$         & $z_0,z_1,z_2$\\
  $G_6$ & 3 & $x_0,x_1$ & $1,\zeta_3,\zeta_3^2$ & $z_0,z_1,z_2$\\
  $G_4$ & 6 & $1,-1$    & $1,\zeta_3,\zeta_3^2$ & $z_0,z_1,z_2$\\
  \hline
\end{tabular}\
$ $\\ $ $\\
$\textrm{\small Specializations of the parameters for }\mathcal{H}(G_7)$
\end{center}

\subsection*{The groups $G_8$, $G_9$, $G_{10}$, $G_{11}$, $G_{12}$, $G_{13}$, $G_{14}$, $G_{15}$}

The following table gives the specializations of the parameters of
the generic Hecke algebra $\mathcal{H}(G_{11})$,
$(x_0,x_1;y_0,y_1,y_2;z_0,z_1,z_2,z_3)$, which give the generic
Hecke algebras of the groups $G_8,\ldots,G_{15}$ (\cite{Ma2}, Table
4.9).

\begin{center}
\begin{tabular}{|c|c|c|c|c|}
  \hline
  Group & Index & S & T & U \\
  \hline
  $G_{11}$ & 1  & $x_0,x_1$ & $y_0,y_1,y_2$         & $z_0,z_1,z_2,z_3$\\
  $G_{10}$ & 2  & $1,-1$    & $y_0,y_1,y_2$         & $z_1,z_1,z_2,z_3$\\
  $G_{15}$ & 2  & $x_0,x_1$ & $y_0,y_1,y_2$         & $\sqrt{u_0},\sqrt{u_1},-\sqrt{u_0},-\sqrt{u_1}$\\
  $G_9$    & 3  & $x_0,x_1$ & $1,\zeta_3,\zeta_3^2$ & $z_0,z_1,z_2,z_3$\\
  $G_{14}$ & 4  & $x_0,x_1$ & $y_0,y_1,y_2$         & $1,i,-1,-i$\\
  $G_8$    & 6  & $1,-1$    & $1,\zeta_3,\zeta_3^2$ & $z_0,z_1,z_2,z_3$\\
  $G_{13}$ & 6  & $x_0,x_1$ & $1,\zeta_3,\zeta_3^2$ & $\sqrt{u_0},\sqrt{u_1},-\sqrt{u_0},-\sqrt{u_1}$\\
  $G_{12}$ & 12 & $x_0,x_1$ & $1,\zeta_3,\zeta_3^2$ & $1,i,-1,-i$\\
  \hline
\end{tabular}\
$ $\\ $ $\\
$\textrm{\small Specializations of the parameters for }\mathcal{H}(G_{11})$
\end{center}

\subsection*{The groups $G_{16}$, $G_{17}$, $G_{18}$, $G_{19}$, $G_{20}$, $G_{21}$, $G_{22}$}

The following table gives the specializations of the parameters of
the generic Hecke algebra $\mathcal{H}(G_{19})$,
$(x_0,x_1;y_0,y_1,y_2;z_0,z_1,z_2,z_3,z_4)$, which give the generic
Hecke algebras of the groups $G_{16},\ldots,G_{22}$ (\cite{Ma2},
Table 4.12).

\begin{center}
\begin{tabular}{|c|c|c|c|c|}
  \hline
  Group & Index & S & T & U \\
  \hline
  $G_{19}$ & 1  & $x_0,x_1$ & $y_0,y_1,y_2$         & $z_0,z_1,z_2,z_3,z_4$\\
  $G_{18}$ & 2  & $1,-1$    & $y_0,y_1,y_2$         & $z_0,z_1,z_2,z_3,z_4$\\
  $G_{17}$ & 3  & $x_0,x_1$ & $1,\zeta_3,\zeta_3^2$ & $z_0,z_1,z_2,z_3,z_4$\\
  $G_{21}$ & 5  & $x_0,x_1$ & $y_0,y_1,y_2$         & $1,\zeta_5,\zeta_5^2,\zeta_5^3,\zeta_5^4$\\
  $G_{16}$ & 6  & $1,-1$    & $1,\zeta_3,\zeta_3^2$ & $z_0,z_1,z_2,z_3,z_4$\\
  $G_{20}$ & 10 & $1,-1$    & $y_0,y_1,y_2$         & $1,\zeta_5,\zeta_5^2,\zeta_5^3,\zeta_5^4$\\
  $G_{22}$ & 15 & $x_0,x_1$ & $1,\zeta_3,\zeta_3^2$ & $1,\zeta_5,\zeta_5^2,\zeta_5^3,\zeta_5^4$\\
  \hline
\end{tabular}\
$ $\\ $ $\\
$\textrm{\small Specializations of the parameters for }\mathcal{H}(G_{19})$
\end{center}

\subsection*{The groups $G_{25}$, $G_{26}$}

The following table gives the specialization of the parameters of
the generic Hecke algebra $\mathcal{H}(G_{26})$,
$(x_0,x_1;y_0,y_1,y_2)$, which give the generic Hecke algebra of the
group $G_{25}$ (\cite{Ma5}, Theorem 6.3).

\begin{center}
\begin{tabular}{|c|c|c|c|}
  \hline
  Group & Index & S & T \\
  \hline
  $G_{26}$ & 1  & $x_0,x_1$ & $y_0,y_1,y_2$\\
  $G_{25}$ & 2  & $1,-1$    & $y_0,y_1,y_2$\\
  \hline
\end{tabular}\
$ $\\ $ $\\
$\textrm{\small Specialization of the parameters for }\mathcal{H}(G_{26})$
\end{center}

\end{document}